\magnification 1200

\font\tenfrak=eufm10
\font\sevenfrak=eufm7
\font\fivefrak=eufm5
\newfam\frakfam
\textfont\frakfam=\tenfrak
\scriptfont\frakfam=\sevenfrak
\scriptscriptfont\frakfam=\fivefrak
\def\frak{\fam\frakfam}
\def\eu{{\frak eu}}
\def\hh{\hskip-.5pt}
\def\sss#1{\langle#1\rangle}
\def\CP{{\frak CP}}
\def\dist{\mathop{\rm dist}\nolimits}
\def\Arg{\mathop{\rm Arg}\nolimits}
\def\id{\mathop{\rm id}\nolimits}

\font\tenams=msbm10
\font\sevenams=msbm7
\font\fiveams=msbm5
\newfam\amsfam
\textfont\amsfam=\tenams
\scriptfont\amsfam=\sevenams
\scriptscriptfont\amsfam=\fiveams
\def\Bbb#1{{\fam\amsfam #1}}
\def\Cal#1{{\cal #1}}

\font\amsfa=msam10
\def\sq{\hfill\hbox{\amsfa\char'003}}

\centerline{\bf Euler characteristic  of the configuration space of a
complex\footnote{$^\dagger$}{\rm Published in Colloquium Mathematicum
Volume 89, Issue 1, 2001}}
\footnote{}{2000 {\it Mathematics Subject Classification :} 57}
\smallskip
\centerline{\'Swiatos\l aw R. Gal}
\font\tt=cmtt10
\centerline{\tt http://www.math.uni.wroc.pl/\~{}sgal/}
\bigskip
{
\narrower\narrower\smallskip\noindent
\sevenbf Abstract. \sevenrm
A closed form formula (generating function)
for the Euler characteristic
of the configuration space of $\scriptstyle n$
particles in a simplicial complex is given.
\smallskip}
\parskip=\smallskipamount\parindent=0pt

\bigskip
{\bf Introduction.}
\medskip

The Euler characteristic, being independent
of the structure of a
simplicial complex on a topological space,
may be computed in terms of zeroth order data, i.e.
the number of cells of each dimension. We prove that
the Euler characteristics of configuration spaces
may be computed in terms of first order data of
the underlying space, i.e.~the Euler characteristics
of links of cells.

When computing the Euler characteristic
of the configuration space on $n$ (ordered)
particles in a simplicial complex $X$
we mimic the inductive approach which gives
the desired formula in the case of manifolds.
Consider the projection
$\Cal C_n(X)\to X$ onto the first particle.
Even though this map is not in general a fibration,
we can relate the Euler characteristics of the
base, the total space and the fibers of this map.
A major feature here is a very interesting measure on $X$.

Next we restate the formula as a differential equation
for the exponential generating function, which is then solved
explicitly.

The result (and overall strategy) generalizes
a similar result for the case of graphs, which was known to
M. W.~Davis, H.~Glover, T.~Januszkiewicz and J.~\'Swi%
{\newdimen\bsp
\bsp=.10em
\advance\bsp by12\fontdimen1\the\font\setbox0=\hbox{a}\setbox1=\hbox 
  to\wd0{a\hss\kern-\bsp\lower1.4ex\hbox{\char96}\hss}\ht1=\ht0
  \dp1=\dp0\box1}tkowski
around 1997.

\bigskip
{\bf 0. Definitions and conventions.}
\medskip

For any topological space $X$ define the
{\it configuration space of\/ $n$ (ordered\/) particles in\/ $X$} as
$\Cal C_n(X) : =X^n-\{\xi: \xi_i=\xi_j \hbox{ for some }i\neq j\}$;
$\xi_i$ will be referred to as  the $i$th {\it particle} of the
{\it configuration $\xi$}.
Define $\chi_n(X) : =\chi \Cal C_n(X)$
to be the Euler characteristic of $\Cal C_n(X)$.

A subspace $S\subset X$ is {\it collared} if there
is a tubular neighborhood $V\simeq S\times (-1,1)$ of $S$ in $X$.
We say that a complex $Y$ is obtained by
{\it C(ut)\&P(aste) surgery\/}
from a complex $X$ if one can find a collared subcomplex $S\subset X$
(one may think that $X$ is a barycentric subdivision of some complex $Z$ and
$S$ is transversal to $Z$, i.e.~when restricted to any cell of $Z$ it has
codimension $1$)
such that $Y$ can be obtained from $X$ by cutting along
$S$ and then gluing again by some cellular automorphism of $S$.
We define the  {\it C\&P equivalence relation\/} on
the category of finite complexes
as the equivalence relation generated by C\&P surgery and subdivision.
By $\CP$ we denote the Grothendieck ring generated by C\&P classes of
finite complexes with
disjoint union as
sum and with Cartesian product as multiplication
(the class of $X$ will be denoted $[X]$).

Similarly,
in the topological category we define t(\hh opological\hh)C\&P surgery,\hskip-.5pt
and the (topological) $t\CP$ ring.

\proclaim Fact. The Euler characteristic is {\it tC\&P} invariant
and may be regarded as a ring homomorphism $\chi \colon t\CP\to\Bbb Z$.

{\it Proof :} From the exact sequence of a pair it follows that
$\chi(X)=\chi(V)-\chi(X,V)$ and from the excision lemma
$\chi(X,V)=\chi(X-S,V-S$). Since the pairs
$(X-S,V-S)$ and $(Y-S,V-S)$ are homeomorphic
the proof is complete.\sq

We will often use a special surgery (called
an {\it amputation}):
let  $A\subset Z$ be an open subset with collared boundary; take
$X=Z\sqcup(I\times \partial A)$,
$S=\partial A\sqcup(\{1/2\}\times\partial A)$
and let the automorphism
interchange both copies of $\partial A$; then
$Y=(Z-A)\sqcup \bar{A}{I\,}$, therefore
$[Z]=[Z-A]+[\bar{A}]-[I\times\partial A]$.

\bigskip
{\bf 1. The integral formula.}
\medskip

From now on,
let $X$ denote a finite simplicial complex.
Let $CX$ denote the cone over $X$.
If $\sigma$ is a cell of $X$, define:

{\parindent=.4in
\item{$\bullet$} $d_\sigma$ to be the dimension of $\sigma$,
\item{$\bullet$} $L_\sigma$ to be the normal link of $\sigma$,
\item{$\bullet$} $\sss\sigma$ to be a
sufficiently small contractible (open) neighborhood
of the center of $\sigma$.

}

\proclaim Theorem 1. Let\/ $\mu$ be the measure$^{1)}$
on\/ $X$ given by\/
$\mu(\sigma)=1-\chi L_\sigma$ for any cell\/ $\sigma$. Then
$$
\chi_n(X)=\int_X \chi_{n-1}(X-\sss\sigma)\mu(\sigma).\leqno(*)
$$

\footnote{}{\llap{$^{1)}$\ }One
should take the $\sigma$-algebra generated by
the cells of $X$; on the other
hand $\mu$ may be understood as a functional on the space of functions
which are constant on each open cell.}

{\it Remark :} The above formula is obvious when $X$ is
a manifold, since then $\pi$
is a bundle and the Euler characteristic of the total space is
the product of the Euler characteristics of the fiber and the base.
In this special case ($*$) reads: $\chi_n(X)=\chi_{n-1}(X-U)\chi(X)$,
where $U$ is a small open ball.

In the general
case the fiber of $\pi$ changes when $\xi_1$ approaches faces
of non-zero codimension.
The idea of the proof is to write a
deformation retraction that corrects
the picture at the vertices, use amputation surgery,
and continue inductively.

{\it Proof of Theorem 1. :}
Consider on $X$ an auxiliary metric $d$ compatible with
the simplicial structure,
such that each cell is a regular simplex with a side of length $8$.
Let $B_\delta(x)$ and $B_\delta(A)$ denote
the open metric balls of radius
$\delta$ around $x\in X$ and $A\subset X$ respectively.

We
will denote by $\pi:{\cal C}_n(X)\to X$ the projection 
onto the first factor.

Let us define a {\it projection} $p_k$ onto the $k$ skeleton $X^{(k)}$
(it is defined only for some $x$ from $X$, namely for those which are
close to $X^{(k)}$, and far from $X^{(k-1)}$). We set $p_k(x)$ to be
the point of $X^{(k)}$ nearest to $x$.
\vadjust{\eject}

Put $\varepsilon : =1/8$;
define inductively sequences of spaces:
\vskip 4pt

(1) $A_{k}\subset \Cal C_n(X)$ consists of the
     configurations such that
     $\dist(\xi_1,X^{(l)})\geq\varepsilon^l$ for $l<k$
     (roughly speaking $\xi_1$ is far from $X^{(k-1)}$).

(2) $D_k\subset A_k$ consists of
the configurations such that
     if $\dist(\xi_1,X^{(k)})<\varepsilon^k$ then $\xi_1$ is the only
     particle in $B_{\varepsilon^k}(p_k(\xi_1))$.
\vskip 4pt

 For any $\sigma$ such that $d_\sigma=k$
     put $N_\sigma=B_{\varepsilon^k/2}(\sigma)\cap\pi(D_k)$
     ($N_\sigma$ is homeomorphic to $\sigma\times CL_\sigma$,
     its boundary in $\pi(D_k)$ is homeomorphic to
     $\sigma\times L_\sigma$).
Then $\pi^{-1}(N_\sigma)\cap {D_k}$ is homeomorphic to the product
     $N_\sigma\times \Cal C_{n-1}(X-\sss\sigma)$.
     Amputating all the $N_\sigma$'s we write
$$
 [D_k]=[E_k]+\sum_{\sigma :\, d_\sigma=k}[\sigma]
     ([CL_\sigma]-[L_\sigma\times I])[\Cal C_{n-1}(X-\sss\sigma)],
$$
 and thus
$$
\chi D_k=\chi E_k+\sum_{\sigma : \,
d_\sigma=k}(1-\chi L_\sigma)
    \chi_{n-1}(X-\sss\sigma),
$$
where $E_k=\pi^{-1}(X-\bigcup N_\sigma)\cap{D_k}$.

We will complete the proof when we show that $\chi D_k=\chi E_{k-1}$.
This will be done in two steps.

By {\it polar coordinates}\/ on
$CY=[0,\delta]\times Y/\{0\}\times Y$ (we will also denote the class of
$\{0\}\times Y$ by $0$) we mean a pair of functions: the
{\it modulus}\/
$|\cdot| :  CY\to[0,t]$, and the
{\it argument}\/ $\Arg : (CY-0)\to Y$
which are the projections onto factors.

{\sl Step} 1.
{\it $D_k$ is a deformation retract of\/
$A_k$ (thus\/ $\chi A_k=\chi D_k$).}

We can assume that $\xi_1\in B_{\varepsilon^k}(\sigma)$.
Note that $B_{\varepsilon^k}(p_k(\xi_1))$ is a cone over its boundary
(in $X$). Consider polar coordinates on it.

For configuration $\xi$ define
$\varrho : =\min\{d(\xi_i,p_k(\xi_1)) :  i>1\}$.
Both $\xi_1$ and some other particle are
in $B_{\varepsilon^k}(p_k(\xi_1))$ iff $\varrho<\varepsilon^k$.
If this is the case (if not,
then $\xi\in D_k$) define $r : =\max(\varrho,d(\xi_1,p_k(\xi_1)))$.

We need a continuous function
$n : [0,2\varepsilon^k]\times(0,\varepsilon^k]\to[0,2\varepsilon^k]$
such that
$n(\cdot,R)$ is a homeomorphism
(preserving endpoints) of $[0,2\varepsilon^k]$
for any $R$, $n(R,R)=\varepsilon^k$, and $n(\cdot,\varepsilon^k)=\id$.
An example of such a function is
$$
n(t,R) : ={2\varepsilon^{2k}-R^2\over R(2\varepsilon^k-R)}t+
{R-\varepsilon^k\over R(2\varepsilon^k-R)}t^2.
$$

Define
$$
\lambda_{\xi,\theta}(y) : =\cases{
n(|y|,(1-\theta)r+\theta\varepsilon^k)\hbox{Arg}(y)&
if $y\in B_{\varepsilon^k}(p_k(\xi_1))$,\cr
y&otherwise.\cr}
$$
For each $\theta$ this is a homeomorphism of
$X$ continuously depending on $\xi$.

The desired deformation retraction is now given
by the formula
$$
{\mit\Lambda}_\theta(\xi_1,\ldots,\xi_n) : =
(\lambda_{\xi,\theta}(\xi_1),\ldots,\lambda_{\xi,\theta}(\xi_n)).
$$

{\sl Step} 2.
{\it $A_{k+1}$ is a deformation retract of\/ $E_k$\
(thus\/ $\chi E_k=\chi A_{k+1}$).}

For simplicity of notation we will assume that
$B_{2\varepsilon^k}(\sigma)$ is a metric product
$(\sigma\cap\pi(D_k))\times CL_\sigma$ (to do this one has to change
the metric a bit in the neighborhood of $\sigma$).
Consider polar coordinates on
$CL_\sigma$.

We need a continuous function
$n : [\varepsilon^k,2\varepsilon^k]\times[0,\varepsilon^k/2)\to
[\varepsilon^k,2\varepsilon^k]$ such that
$n(\cdot,R)$ is an embedding of
$[\varepsilon^k,2\varepsilon^k]$ into itself,
$n(2\varepsilon^k,\cdot)=2\varepsilon^k$, if $R>0$ then
$n(\cdot,R)>\varepsilon^k$, and $n(\cdot,0)=\id$.
An example of such a function is
$$
n(t,R) : ={\varepsilon^{k}-R\over\varepsilon^k}\,t+
{2R\over\varepsilon^k}.
$$

Define
$$
\lambda_{\xi,\theta}(y) : =\cases{
(p_k(y),n(|y|,\theta(\varepsilon^k-|\xi_1|))\hbox{Arg}(y))&
if $y\in B_{2\varepsilon^k}(\sigma)$,\cr
y&otherwise,\cr}
$$
and
$$
\zeta_\theta:=
(p_k(\xi_1),(|\xi_1|+\theta(\varepsilon^k-|\xi_1|))\hbox{Arg}(\xi_1)).
$$

The desired deformation retraction is now given
by the formula
$$
{\mit\Lambda}_\theta(\xi_1,\ldots,\xi_n) : =(\zeta_\theta,
\lambda_{\xi,\theta}(\xi_2),\ldots,\lambda_{\xi,\theta}(\xi_n)).\sq
$$

\bigskip
{\bf 2. Generating function}
\medskip

\proclaim Definition.
Let $\eu_X(t) : =\sum \chi_n(X)t^n/n!$.

We use this particular generating function since it
satisfies an interesting differential equation. Namely, it is clear that
Theorem 1 can be restated as the differential equation
$$
\eu_X'=\int_X\eu_{X-\sss\sigma}\mu(\sigma).\leqno(\star)
$$

A preliminary step to solve $(\star)$ is

\proclaim Proposition 1. $\eu$ is a group homomorphism of the
additive group of\/ $\CP$
to the group of units\/ $\Bbb Z[[t]]^\times$.

{\it Proof :}
First we prove that this is a homomorphism.
For any configuration of $n$ particles in
$X\sqcup Y$ there is a partition
$I\sqcup J=\{1,\dots,n\}$ such that $x_i\in X$ for $i\in I$
and $x_j\in Y$ for $j\in J$, i.e.
$$
\Cal C_n(X\sqcup Y)\simeq\bigsqcup_{I\sqcup J
=\{1,\dots,n\}}\Cal C_{\#I}(X)
\Cal C_{\#J}(Y).
$$
Therefore
$$
{\chi_n(X\sqcup Y)\over n!}=
{1\over n!}\sum_k { n\choose  k}
 \chi_k(X)\chi_{n-k}(Y) =
\sum_k
{\chi_k(X)\over k!}\cdot
{\chi_{n-k}(Y)\over (n-k)!}.
$$

We prove inductively that $\chi_n$ is C\&P invariant.
Let $Y$ be obtained from $X$ by C\&P surgery along $S$. To show
that $\chi_n(X)=\chi_n(Y)$ using $(\star)$,
it is sufficient to check that
for any simplex $\sigma$ its measure in $X$ and $Y$ is the same,
and $X-\sss\sigma$ and $Y-\sss\sigma$
are C\&P related. This is obvious when $\sigma$ is not in $S$.
If $\sigma$ does belong to $S$ then its link in $X$ is the same
as a link in $Y$
(it is the suspension of a link in $S$),
the second statement follows from the fact
that a neighborhood of $S$ is isomorphic to $S\times(0,1)$, and
$\sss\sigma$
can be pushed away from $S$ by some automorphism of $S\times (0,1)$ which
is the identity near the ends of the interval. Then we can apply C\&P
surgery along $S$ and eventually move $\sss\sigma$ back.\sq

\bigskip
{\bf 3. Explicit formula for $\eu$.}
\medskip

\proclaim Proposition 2.
$$
{\eu_{CX}(t)\over\eu_{X\times I}(t)}=1+(1-\chi X)t.
$$

{\it Remark :} The results of
the previous sections are still valid if we
assume that cells of $X$ are products of simplices.
Therefore we can put a product structure on $X\times I$. The cells of
$CX$ are cones of cells of $X$ and their sides.

{\it Proof of Proposition 2. :}
Let $\varrho$ be a cell of $CX$ different from the vertex of the cone.
Then
$$
[CX-\sss\varrho]=[CX]+[X\times I-\sss\varrho]-[X\times I],
$$
i.e.
$$
\eu_{CX-\sss\varrho}={\eu_{CX}\over\eu_{X\times I}}
\eu_{X\times I-\sss\varrho}.
$$
If $\varrho$ is a vertex of a cone then $CX-\sss\varrho=X\times I$.
Combining the above equalities with ($\star$) we obtain
$$
\eqalign{
\eu_{CX}'
&=
(1-\chi X)\eu_{X\times I}
+\int_{X\times[0,1)}{\eu_{CX}\over\eu_{X\times I}}
\eu_{X\times I-\sss\varrho}\mu(\varrho)\cr
&=
(1-\chi X)\eu_{X\times I}+{\eu_{CX}\eu_{X\times I}'\over\eu_{X\times I}}.
\cr}
$$
The last equality follows since
$\mu(\sigma\times\{1\})=0$ (the link is a cone)
and the integral may be taken over $X\times I$.
Eventually we obtain
$$
\bigg({\eu_{CX}\over\eu_{X\times I}}\bigg)'
=1-\chi X.\sq
$$

\proclaim Theorem 2. Let\/ $X$ be a complex. For any cell\/ $\sigma$ let\/
$d_\sigma$ and\/ $v_\sigma$
denote  respectively the dimension of\/
$\sigma$ and the Euler characteristic of\/
the normal link\/ $L_\sigma$ of\/ $\sigma$. Then
$$
\eu_X(t)=\prod_{\sigma}
(1+(-1)^{d_\sigma}(1-v_\sigma)t )^{(-1)^{d_\sigma}}.\leqno(1)
$$

{\it Proof :}
Apply $\eu$ to the equality
$[X-\sss\sigma]=[X]+[\partial \sss\sigma\times I]-[\overline{\sss\sigma}]$,
and integrate it over $X$ to obtain
$$
\eu_X'=\eu_X
\int_X{\eu_{\partial \sss\sigma\times I}\over\eu_{\overline{\sss\sigma}}}
\mu(\sigma).
$$
Since $\overline{\sss\sigma}=C\partial \sss\sigma$, we apply Proposition 2:
$$\leqalignno{
(\log\eu_X)'(t)
&=
\int_X{1\over1+(1-\chi \partial \sss\sigma)t}\mu(\sigma)&(2)\cr
&=
\int_X{{d\over dt}\log(1+(1-\chi \partial \sss\sigma)t)\over
1-\chi \partial \sss\sigma}\mu(\sigma).
\cr}
$$
Also $\partial \sss\sigma$ is the
$d_\sigma$-fold suspension of
$L_\sigma$, so
$1-\chi \partial \sss\sigma=(-1)^{d_\sigma}(1-v_\sigma)$.
Therefore, rewriting the right hand side we have
$$
(\log\eu_X)'(t)=\sum_{\sigma}
(-1)^{d_\sigma}{d\over dt}\log(1+(-1)^{d_\sigma}(1-v_\sigma)t),
$$
which proves the theorem.\sq

{\it Note :} When $\chi\partial \sss\sigma=1$,
the corresponding summand in
$(2)$ is $0$ and the factor in $(1)$ is $1$.

\bigskip
{\bf Acknowledgements.}
\medskip

I am grateful to my advisor Tadeusz Janusz\-kie\-wicz
for suggesting the project and help throughout.
I thank Jan Dymara for carefully reading the paper
and for his help in improving the presentation.

\bigskip
\font\sc=cmcsc10\sc
Institute of Mathematics, University of Wroc\l aw,
Pl. Grunwaldzki 2/4, 50-384 Wroc\l aw, Poland,
e-mail : \tt sgal@math.uni.wroc.pl

\bye